\documentclass[a4paper, 12pt, leqno]{article}
\usepackage{amsmath}
\usepackage{amscd}
\usepackage{amssymb}
\usepackage{amsthm}
\newtheorem{lemma}{Lemma}[section]
\newtheorem{corollary}[lemma]{Corollary}
\newtheorem{proposition}[lemma]{Proposition}
\newtheorem{theorem}[lemma]{Theorem}
\newtheorem{remark}[lemma]{Remark}

\def\Gb{{\mathfrak{b}}}
\def\Gg{{\mathfrak{g}}}

\def\BA{{\mathbb{A}}}
\def\BC{{\mathbb{C}}}
\def\BD{{\mathbb{D}}}

\def\BN{{\mathbb{N}}}

\def\BQ{{\mathbb{Q}}}

\def\BZ{{\mathbb{Z}}}

\def\lan{\langle}
\def\ran{\rangle}

\def\pt{{\rm{pt}}}

\def\codim{\mathop{\rm codim}\nolimits}

\def\Hom{\mathop{\rm Hom}\nolimits}
\def\id{\mathop{\rm id}\nolimits}
\def\Id{\mathop{\rm Id}\nolimits}

\def\op{{\mathop{\rm op}\nolimits}}
\def\Ob{\mathop{\rm Ob}\nolimits}

\def\Spec{\mathop{\rm Spec}\nolimits}

\newcommand{\iso}{\mathrel{\longrightarrow{\kern-18pt\raise3.5pt\hbox{$\sim$
}}}\enspace}
\def\X{{X'}}
\def\Y{{Y'}}
\def\HS{\mathop{\rm HS}\nolimits}
\def\HM{\mathop{\rm HM_d}\nolimits}
\def\tHM{\mathop{\rm HM_c}\nolimits}
\def\HMB{\mathop{{\rm HM}^B_{\rm d}}\nolimits}
\def\tHMB{\mathop{{\rm HM}^{B^-}_{\rm c}}\nolimits}
\def\BDD{\mathop{\BD_{\rm d}}\nolimits}
\def\BDC{\mathop{\BD_{\rm c}}\nolimits}

\title{Parabolic Kazhdan-Lusztig polynomials\\ and Schubert varieties}
\author{Masaki Kashiwara%
\thanks{Research Institute for Mathematical Sciences, Kyoto University,
Kyoto, 606--8502, Japan} and
Toshiyuki Tanisaki%
\thanks{Department of Mathematics, Hiroshima University,
Higashi-Hiroshima, 739--8526, Japan}}
\begin{document}

\maketitle
\begin{abstract}
We shall give a description of the intersection cohomology groups of the
Schubert varieties in partial flag manifolds over symmetrizable Kac-Moody
Lie algebras in
terms of parabolic Kazhdan-Lusztig polynomials introduced by Deodhar.
\end{abstract}

\section{Introduction}
\setcounter{equation}{0}
For a Coxeter system $(W,S)$ Kazhdan-Lusztig~\cite{KL1}, \cite{KL2}
introduced polynomials
\[
P_{y,w}(q)=\sum_{k\in\BZ}P_{y,w,k}q^k\in\BZ[q],\qquad
Q_{y,w}(q)=\sum_{k\in\BZ}Q_{y,w,k}q^k\in\BZ[q],
\]
called a Kazhdan-Lusztig polynomial and an inverse  Kazhdan-Lusztig
polynomial respectively.
Here, $(y,w)$ is a pair of elements of $W$ such that $y\leqq w$ with respect
to the Bruhat order.
These polynomials play important roles in various aspects of the
representation theory of reductive algebraic groups.

In the case $W$ is associated to a symmetrizable Kac-Moody Lie algebra
$\Gg$, the polynomials have the following geometric meanings.
Let $X=G/B$ be the corresponding flag variety (see Kashiwara~\cite{K1}), and
set $X^w=B^-wB/B$ and $X_w=BwB/B$ for $w\in W$.
Here $B$ and $B^-$ are the ``Borel subgroups" corresponding to the standard
Borel subalgebra $\Gb$ and its opposite $\Gb^-$ respectively.
Then $X^w$ (resp.\ $X_w$) is an $\ell(w)$-codimensional (resp.\
$\ell(w)$-dimensional) locally closed subscheme of the infinite-dimensional
scheme $X$.
Here $\ell(w)$ denotes the length of $w$ as an element of the Coxeter group
$W$.
Set $X'=\bigcup_{w\in W}X_w$.
Then $X'$ coincides with the flag variety considered by
Kac-Peterson~{\cite{KP}}, Tits~{\cite{T}}, et al.
Moreover we have
\[
X=\bigsqcup_{w\in W}X^w,\qquad X'=\bigsqcup_{w\in W}X_w,
\]
and
\[
\overline{X^w}=\bigsqcup_{y\geqq w}X^y, \qquad
\overline{X_w}=\bigsqcup_{y\leqq w}X_y
\]
for any $w\in W$.

By Kazhdan-Lusztig~\cite{KL2} we have the following result (see also
Kashiwara-Tanisaki~\cite{KT1}).
\begin{theorem}
\label{thm:intro:1}
\begin{itemize}
\item[\rm(i)]
Let $w, y\in W$ satisfying $w\leqq y$.
Then we have
\[
H^{2k+1}({}^\pi\BQ^H_{X^w})_{yB/B}=0,\qquad
H^{2k}({}^\pi\BQ^H_{X^w})_{yB/B}=\BQ^H(-k)^{\oplus Q_{w,y,k}}
\]
for any $k\in\BZ$.
\item[\rm(ii)]
The multiplicity of the irreducible Hodge module
${}^\pi\BQ^H_{X^y}[-\ell(y)](-k)$ in the
Jordan H\"older series of the Hodge module $\BQ^H_{X^w}[-\ell(w)]$ coincides
with $P_{w,y,k}$.
\end{itemize}
\end{theorem}
\begin{theorem}
\label{thm:intro:2}
\begin{itemize}
\item[\rm(i)]
Let $w, y\in W$ satisfying $w\geqq y$.
Then we have
\[
H^{2k+1}({}^\pi\BQ^H_{X_w})_{yB/B}=0,\qquad
H^{2k}({}^\pi\BQ^H_{X_w})_{yB/B}=\BQ^H(-k)^{\oplus P_{y,w,k}}
\]
for any $k\in\BZ$.
\item[\rm(ii)]
The multiplicity of the irreducible Hodge module
${}^\pi\BQ^H_{X_y}[\ell(y)](-k)$ in the
Jordan H\"older series of the Hodge module $\BQ^H_{X_w}[\ell(w)]$ coincides
with $Q_{y,w,k}$.
\end{itemize}
\end{theorem}
Here ${}^\pi\BQ^H_{X^w}[-\ell(w)]$ and ${}^\pi\BQ^H_{X_w}[\ell(w)]$ denote
the Hodge modules corresponding to the perverse sheaves
${}^\pi\BQ_{X^w}[-\ell(w)]$ and ${}^\pi\BQ_{X_w}[\ell(w)]$ respectively.
In Theorem~\ref{thm:intro:1} we have used the convention so that
${}^\pi\BQ^H_{Z}[-\codim Z]$ is a Hodge module for a locally closed
finite-codimensional subvariety $Z$ since we deal with sheaves supported on
finite-codimensional subvarieties, while in Theorem~\ref{thm:intro:2} we
have used another convention so that  ${}^\pi\BQ^H_{Z}[\dim Z]$ is a Hodge
modules for a locally closed finite-dimensional subvariety $Z$ since we deal
with sheaves supported on finite-dimensional subvarieties.

Let $J$ be a subset of $S$.
Set $W_J=\langle J\rangle$ and denote by $W^J$ the set of elements $w\in W$
whose length is minimal in the coset $wW_J$.
In \cite{D} Deodhar introduced two generalizations of the Kazhdan-Lusztig
polynomials to this relative situation.
For $(y,w)\in W^J\times W^J$ such that $y\leqq w$ we denote the parabolic
Kazhdan-Lusztig polynomial for $u=-1$ by
\[
P^{J,q}_{y,w}(q)=\sum_{k\in\BZ}P^{J,q}_{y,w,k}q^k\in\BZ[q],
\]
and that for $u=q$ by
\[
P^{J,-1}_{y,w}(q)=\sum_{k\in\BZ}P^{J,-1}_{y,w,k}q^k\in\BZ[q]
\]
contrary to the original reference \cite{D}.
We can also define inverse parabolic Kazhdan-Lusztig polynomials
\[
Q^{J,q}_{y,w}(q)=\sum_{k\in\BZ}Q^{J,q}_{y,w,k}q^k\in\BZ[q],
\qquad
Q^{J,-1}_{y,w}(q)=\sum_{k\in\BZ}Q^{J,-1}_{y,w,k}q^k\in\BZ[q]
\]
(see \S~\ref{sec:KL} below)

The aim of this paper is to extend Theorem~\ref{thm:intro:1} and
Theorem~\ref{thm:intro:2} to this relative situation using the partial flag
variety corresponding to $J$.

Let $Y$ be the partial flag variety corresponding to $J$.
Let $1_Y$ be the origin of $Y$ and set $Y^w=B^-w1_Y$ and $Y_w=Bw1_Y$ for
$w\in W^J$.
Then $Y^w$ (resp.\ $Y_w$) is an $\ell(w)$-codimensional (resp.\
$\ell(w)$-dimensional) locally closed subscheme of the infinite-dimensional
scheme $Y$.
Set $Y'=\bigcup_{w\in W^J}Y_w$.
Then we have
\[
Y=\bigsqcup_{w\in W^J}Y^w,\qquad Y'=\bigsqcup_{w\in W^J}Y_w,
\]
and
\[
\overline{Y^w}=\bigsqcup_{y\geqq w}Y^y, \qquad
\overline{Y_w}=\bigsqcup_{y\leqq w}Y_y
\]
for any $w\in W^J$.

We note that the construction of the partial flag variety similar to the ordinary flag variety in Kashiwara~\cite{K1} has not yet appeared in the literature.
In the case where $W_J$ is a finite group (especially when $W$ is an affine Weyl group), we can construct the partial flag variety $Y=G/P$ and the properties of Schubert varieties in $Y$ stated above are established in exactly the same manner as in Kashiwara~\cite{K1} and Kashiwara-Tanisaki~\cite{KT2}.
In the case $W_J$ is an infinite group we can not define the ``parabolic subgroup" $P$ corresponding to $J$ as a group scheme and hence the arguments in Kashiwara~\cite{K1} are not directly generalized.
We leave the necesary modification in the case $W_J$ is an infinite group to the future work.

Our main result is the following.
\begin{theorem}
\label{thm:intro:3}
\begin{itemize}
\item[\rm(i)]
Let $w, y\in W^J$ satisfying $w\leqq y$.
Then we have
\[
H^{2k+1}({}^\pi\BQ^H_{Y^w})_{y1_Y}=0,\qquad
H^{2k}({}^\pi\BQ^H_{Y^w})_{y1_Y}=\BQ^H(-k)^{\oplus Q^{J,-1}_{w,y,k}}
\]
for any $k\in\BZ$.
\item[\rm(ii)]
The multiplicity of the irreducible Hodge module
${}^\pi\BQ^H_{Y^y}[-\ell(y)](-k)$ in the
Jordan H\"older series of the Hodge module $\BQ^H_{Y^w}[-\ell(w)]$ coincides
with $P^{J,-1}_{w,y,k}$.
\end{itemize}
\end{theorem}
\begin{theorem}
\label{thm:intro:4}
\begin{itemize}
\item[\rm(i)]
Let $w, y\in W^J$ satisfying $w\geqq y$.
Then we have
\[
H^{2k+1}({}^\pi\BQ^H_{Y_w})_{y1_Y}=0,\qquad
H^{2k}({}^\pi\BQ^H_{Y_w})_{y1_Y}=\BQ^H(-k)^{\oplus P^{J,q}_{y,w,k}}
\]
for any $k\in\BZ$.
\item[\rm(ii)]
The multiplicity of the irreducible Hodge module
${}^\pi\BQ^H_{Y_y}[\ell(y)](-k)$ in the
Jordan H\"older series of the Hodge module $\BQ^H_{Y_w}[\ell(w)]$ coincides
with $Q^{J,-1}_{y,w,k}$.
\end{itemize}
\end{theorem}
In Theorem~\ref{thm:intro:3} we have used the convention so that
${}^\pi\BQ^H_{Z}[-\codim Z]$ is a Hodge module for a locally closed
finite-codimensional subvariety $Z$, and in Theorem~\ref{thm:intro:4} we
have used another convention so that  ${}^\pi\BQ^H_{Z}[\dim Z]$ is a Hodge
modules for a locally closed finite-dimensional subvariety $Z$.

We note that a result closely related to Theorem~\ref{thm:intro:4} was
already obtained by Deodhar~\cite{D}.

The above results imply that the coefficients of the four (oridnary or
inverse) parabolic Kazhdan-Lusztig polynomials are all non-negative in the
case $W$ is the Weyl group of a symmetrizable Kac-Moody Lie algebra.

We would like to thank B.\ Leclerc for leading our attention  to this problem.
We also thank H.\ Tagawa for some helpful comments on the manuscript.
\section{Kazhdan-Lusztig polynomials}
\label{sec:KL}
\setcounter{equation}{0}
Let $R$ be a commutative ring containing $\BZ[q,q^{-1}]$ equipped with a
direct sum decomposition  $R=\bigoplus_{k\in\BZ}R_k$ into $\BZ$-submodules
and an involutive ring endomorphism $R\ni r\mapsto \overline{r}\in R$
satisfying the following conditions:
\begin{equation}
\label{eq:cond:R}
R_iR_j\subset R_{i+j},\quad
\overline{R_i}=R_{-i},\quad
1\in R_0, \quad
q\in R_2,\quad
\overline{q}=q^{-1}.
\end{equation}

Let $(W,S)$ be a Coxeter system.
We denote by $\ell:W\to\BZ_{\geqq0}$ and $\geqq$ the length function and the
Bruhat order respectively.
The Hecke algebra $H=H(W)$ over $R$ is an $R$-algebra with free $R$-basis
$\{T_w\}_{w\in W}$ whose multiplication is determined by the following:
\begin{align}
\label{eq:Hecke:1}
&T_{w_1}T_{w_2}=T_{w_1w_2}\qquad
\mbox{if $\ell(w_1w_2)=\ell(w_1)+\ell(w_2)$},
\\
\label{eq:Hecke:2}
&(T_s+1)(T_s-q)=0 \qquad
\mbox{for $s\in S$.}
\end{align}
Note that $T_e=1$ by \eqref{eq:Hecke:1}.

We define involutive ring endomorphisms
$H\ni h\mapsto\overline{h}\in H$
and $j:H\to H$ by
\begin{equation}
\overline{\sum_{w\in W}r_wT_w}
=\sum_{w\in W}\overline{r}_wT_{w^{-1}}^{-1},
\qquad
j(\sum_{w\in W}r_wT_w)
=\sum_{w\in W}{r}_w(-q)^{\ell(w)}T_{w^{-1}}^{-1}.
\end{equation}
Note that $j$ is an endomorphism of an $R$-algebra.

\begin{proposition}
[Kazhdan-Lusztig~\cite{KL1}]
\label{prop:KL}
For any $w\in W$ there exists a unique $C_w\in H$ satisfying the following
conditions:
\begin{align}
&\mbox{$C_w=\sum_{y\leqq w}P_{y,w}T_y$ with $P_{w,w}=1$ and $P_{y,w}\in
\bigoplus_{i=0}^{\ell(w)-\ell(y)-1}R_i$}\\
&\mbox{for $y<w$,}\nonumber\\
&\overline{C}_w=q^{-\ell(w)}C_w.\label{2.6}
\end{align}
Moreover we have $P_{y,w}\in\BZ[q]$ for any $y\leqq w$.
\end{proposition}
Note that $\{C_w\}_{w\in W}$ is a basis of the $R$-module $H$.
The polynomials $P_{y,w}$ for $y\leqq w$ are called Kazhdan-Lusztig
polynomials.
We write
\begin{equation}
\label{eq:coef:P}
P_{y,w}=\sum_{k\in\BZ}P_{y,w,k}q^k.
\end{equation}

Set $H^*=H^*(W)=\Hom_R(H,R)$.
We denote by $\lan\ ,\ \ran$ the coupling between $H^*$ and $H$.
We define involutions $H^*\ni m\mapsto\overline{m}\in H^*$ and $j:H^*\to
H^*$ by
\begin{equation}\label{2.8}
\langle \overline{m},h\rangle=\overline{\langle m, \overline{h}\rangle},\quad
\langle j(m),h\rangle=\langle m, j(h)\rangle\qquad
\mbox{for $m\in H^*$ and $h\in H$}.
\end{equation}
Note that $j$ is an endomorphism of an $R$-module.
For $w\in W$ we define elements $S_w, D_w\in H^*$  by
\begin{equation}
\langle S_w, T_x\rangle=(-1)^{\ell(w)}\delta_{w,x},
\qquad\langle D_w, C_x\rangle=(-1)^{\ell(w)}\delta_{w,x}.
\end{equation}
Then any element of $H^*$ is uniquely written as an infinite sum
in two ways
$\sum_{w\in W}r_wS_w$
and $\sum_{w\in W}r'_wD_w$ with $r_w$,$r_w'\in R$.
Note that we have
\begin{equation}
S_w=\sum_{y\geqq w}(-1)^{\ell(w)-\ell(y)}P_{w,y}D_y
\end{equation}
by $C_w=\sum_{y\leqq w}P_{y,w}T_y$.
By (\ref{2.6}), we have
\begin{equation}
\overline{D}_w=q^{\ell(w)}D_w,
\end{equation}
and we can write
\begin{equation}
\label{eq:D=sum of S}
D_w=\sum_{y\geqq w}Q_{w,y}S_y,
\end{equation}
where $Q_{w,y}$ are  determined by
\begin{equation}
\label{eq:inv KL}
\sum_{w\leqq y\leqq z}(-1)^{\ell(y)-\ell(w)}Q_{w,y}P_{y,z}=\delta_{w,z}.
\end{equation}
Note that \eqref{eq:D=sum of S} is equivalent to
\begin{equation}
\label{eq:T=sum of C}
T_w=\sum_{y\leqq w}(-1)^{\ell(w)-\ell(y)}Q_{y,w}C_y.
\end{equation}
By \eqref{eq:inv KL} we see easily that
\begin{align}
\label{eq:Q1}
&Q_{w,y}\in\BZ[q],\\
\label{eq:Q2}
&\mbox{$Q_{w,w}=1$ and $\deg Q_{w,y}\leqq{(\ell(y)-\ell(w)-1)/2}$ for $w<y$}.
\end{align}
The polynomials $Q_{w,y}$ for $w\leqq y$ are called inverse Kazhdan-Lusztig
polynomials (see Kazhdan-Lusztig~\cite{KL2}).
We write
\begin{equation}
\label{eq:coef:Q}
Q_{w,y}=\sum_{k\in\BZ}Q_{w,y,k}q^k.
\end{equation}

The following result is proved similarly to Proposition~\ref{prop:KL} (see
Kashiwara-Tanisaki~\cite{KT1}).
\begin{proposition}
\label{prop:inv-KL}
Let $w\in W$.
Assume that $D\in H^*$ satisfies the following conditions:
\begin{align}
&\mbox{$D=\sum_{y\geqq w}r_yS_y$ with $r_w=1$ and $r_y\in
\bigoplus_{i=0}^{\ell(y)-\ell(w)-1}R_i$}\\
&\mbox{for $w<y$,}\nonumber\\
&\overline{D}=q^{\ell(w)}D.
\end{align}
Then we have $D=D_w$.
\end{proposition}

We fix a subset $J$ of $S$ and set
\begin{equation}
W_J=\langle J\rangle,\qquad
W^J=\{w\in W\,;\,ws>w \quad\mbox{for any $s\in J$}\}.
\end{equation}
Then we have
\begin{align}
&W=\bigsqcup_{w\in W^J}wW_J,\\
&\mbox{$\ell(wx)=\ell(w)+\ell(x)$ for any $w\in W^J$ and $x\in W_J$.}
\end{align}
When $W_J$ is a finite group, we denote the longest element of $W_J$ by $w_J$.

Let $a\in\{q,-1\}$ and define $a^\dagger\in\{q, -1\}$ by $aa^\dagger=-q$.
Define an algebra homomorphism $\chi^a:H(W_J)\to R$ by
$\chi^a(T_w)=a^{\ell(w)}$, and denote the corresponding one-dimensional
$H(W_J)$-module by $R^a=R1^a$.
We define the induced module $H^{J,a}$ by
\begin{equation}
H^{J,a}=H\otimes_{H(W_J)}R^a,
\end{equation}
and define $\varphi^{J,a}:H\to H^{J,a}$  by $\varphi^{J,a}(h)=h\otimes1^a$.

It is easily checked that $H^{J,a}\ni k\mapsto \overline{k}\in H^{J,a}$ and
$j^a:H^{J,a}\to H^{J, a^\dagger}$ are well defined by
\begin{equation}
\label{eq:phi,j}
\overline{\varphi^{J,a}(h)}=\varphi^{J,a}(\overline{h}),\quad
j^a(\varphi^{J,a}(h))=\varphi^{J,a^\dagger}(j(h))\qquad
\mbox{for $h\in H$}.
\end{equation}
Note that $j^a$ is a homomorphism of $R$-modules and that
\begin{align}
&\mbox{$\overline{rk}=\overline{r}\overline{k}$\qquad for $r\in R$ and $k\in
H^{J,a}$,}\\
&\mbox{$\overline{\overline{k}}=k$\qquad for $k\in H^{J,a}$,}\\
&j^{a^\dagger}\circ j^{a}=\id_{H^{J,a}}.
\end{align}

For $w\in W^J$ set $T_w^{J,a}=\varphi^{J,a}(T_w)$.
It is easily seen that $H^{J,a}$ is a free $R$-module with basis
$\{T_w^{J,a}\}_{w\in W^J}$.
Note that we have
\begin{equation}
\label{eq:phi-T}
\varphi^{J,a}(T_{wx})=a^{\ell(x)}T^{J,a}_w\qquad
\mbox{for $w\in W^J$ and $x\in W_J$}.
\end{equation}
\begin{proposition}
[Deodhar~\cite{D}]
\label{prop:rel-KL}
For any $w\in W^J$ there exists a unique $C^{J,a}_w\in H^{J,a}$ satisfying
the following conditions.
\begin{align}
&\mbox{$C^{J,a}_w=\sum_{y\leqq w}P^{J,a}_{y,w}T_y$ with $P^{J,a}_{w,w}=1$
and $P^{J,a}_{y,w}\in \bigoplus_{i=0}^{\ell(w)-\ell(y)-1}R_i$}\\
&\mbox{for $y<w$.}\nonumber\\
&\overline{C^{J,a}_w}=q^{-\ell(w)}C^{J,a}_w.
\end{align}
Moreover we have $P^{J,a}_{y,w}\in\BZ[q]$ for any $y\leqq w$.
\end{proposition}
The polynomials $P^{J,a}_{y,w}$ for $y, w\in W^J$ with $y\leqq w$ are called
parabolic Kazhdan-Lusztig polynomials.
We write
\begin{equation}
\label{eq:coef:PJa}
P^{J,a}_{y,w}=\sum_{k\in\BZ}P^{J,a}_{y,w,k}q^k.
\end{equation}

\begin{remark}
{\rm
In the original reference \cite{D} Deodhar uses
\[
(-1)^{\ell(w)}j^{a^\dagger}(C^{J,a^\dagger}_w)=
\sum_{y\leqq w}(-q)^{\ell(w)-\ell(y)}\overline{P^{J,a^\dagger}_{y,w}}T^{J,a}_y
\]
instead of $C^{J,a}_w$ to define the parabolic Kazhdan-Lusztig polynomials.
Hence our $P^{J,a}_{y,w}$ is actually the parabolic Kazhdan-Lusztig
polynomial $P^{J}_{y,w}$ for $u=a^\dagger$ in the terminology of \cite{D}.}
\end{remark}
\begin{proposition}
[Deodhar~\cite{D}]
\label{comparison-P}
Let $w, y\in W^J$ such that $w\geqq y$.
\begin{itemize}
\item[\rm(i)]
We have
\[
P^{J,-1}_{y,w}=\sum_{x\in W_J, yx\leqq w}(-1)^{\ell(x)}P_{yx,w}.
\]
\item[\rm(ii)]
If $W_J$ is a finite group, then we have $P^{J,q}_{y,w}=P_{yw_J,ww_J}$.
\end{itemize}
\end{proposition}
Set
\begin{equation}
H^{J,a,*}=\Hom_R(H^{J,a},R),
\end{equation}
and define ${}^t\varphi^{J,a}:H^{J,a,*}\to H^*$ by
\[
\langle{}^t\varphi^{J,a}(n),h\rangle=\langle n,\varphi^{J,a}(h)\rangle
\qquad
\mbox{for $n\in H^{J,a,*}$ and $h\in H$}.
\]
Then ${}^t\varphi^{J,a}$ is an injective homomorphism of $R$-modules.
We define an involution $-$ of $H^{J,a,*}$ similarly to (\ref{2.8}).
We can easily check that
\begin{equation}
\label{eq:t-phi-bar}
\overline{{}^t\varphi^{J,a}(n)}={}^t\varphi^{J,a}(\overline{n})\qquad
\mbox{for any $n\in H^{J,a,*}$.}
\end{equation}

For $w\in W^J$ we define $S_w^{J,a}, D_w^{J,a}\in H^{J,a,*}$ by
\begin{equation}
\langle S^{J,a}_w, T^{J,a}_x\rangle=(-1)^{\ell(w)}\delta_{w,x},
\qquad\langle D^{J,a}_w, C^{J,a}_x\rangle=(-1)^{\ell(w)}\delta_{w,x}.
\end{equation}
Then any element of $H^{J,a,*}$ is written uniquely as an infinite sum
in two ways
$\sum_{w\in W^J}r_wS^{J,a}_w$ and
$\sum_{w\in W^J}r'_wD^{J,a}_w$
 with $r_w$,$r'_w\in R$.
Note that we have
\begin{equation}
\label{eq:SJa-sum-b}
S^{J,a}_w=\sum_{y\in W^J, y\geqq w}(-1)^{\ell(w)-\ell(y)}P^{J,a}_{w,y}D^{J,a}_y
\end{equation}
by $C^{J,a}_w=\sum_{y\leqq w}P^{J,a}_{y,w}T_y$.
We see easily by (\ref{eq:phi-T}) that
\begin{equation}
\label{eq:SJa-sum}
{}^t\varphi^{J,a}(S^{J,a}_w)=\sum_{x\in W_J}(-a)^{\ell(x)}S_{wx}\qquad
\mbox{for $w\in W^J$.}
\end{equation}
By the definition we have
\begin{equation}
\label{eq:DJa-duality}
\overline{D^{J,a}_w}=q^{\ell(w)}D^{J,a}_w,
\end{equation}
and we can write
\begin{equation}
\label{eq:DJa-sum}
D^{J,a}_w=\sum_{y\in W_J,\,y\geqq w}Q^{J,a}_{w,y}S^{J,a}_y
\end{equation}
where $Q^{J,a}_{w,y}\in R$ are  determined by
\begin{align}
\label{eq:inv-rel-KL}
\sum_{y\in W^J, w\leqq y\leqq z}(-1)^{\ell(y)-\ell(w)}
&Q^{J,a}_{w,y}P^{J,a}_{y,z}=\delta_{w,z}\\
&\mbox{for $w, z\in W^J$ satisfying $w\leqq z$.}
\nonumber
\end{align}
Note that \eqref{eq:DJa-sum} is equivalent to
\begin{equation}
\label{eq:TJa-sum}
T^{J,a}_w=\sum_{y\in W^J, y\leqq
w}(-1)^{\ell(w)-\ell(y)}Q^{J,a}_{y,w}C^{J,a}_y.
\end{equation}
By \eqref{eq:inv-rel-KL} we have for $w,y\in W_J$
\begin{align}
\label{eq:QJa1}
&Q^{J,a}_{w,y}\in\BZ[q],\\
\label{eq:QJa2}
&\mbox{$Q^{J,a}_{w,w}=1$ and $\deg
Q^{J,a}_{w,y}\leqq{(\ell(y)-\ell(w)-1)/2}$ for $w<y$}.
\end{align}
We call the polynomials $Q^{J,a}_{w,y}$ for $w\leqq y$ inverse parabolic
Kazhdan-Lusztig polynomials.
We write
\begin{equation}
\label{eq:coef:QJa}
Q^{J,a}_{w,y}=\sum_{k\in\BZ}Q^{J,a}_{w,y,k}q^k.
\end{equation}

Similarly to Propositions~\ref{prop:KL}, \ref{prop:inv-KL},
\ref{prop:rel-KL}, we can prove the following.
\begin{proposition}
\label{prop:inv-rel-KL}
Let $w\in W^J$.
Assume that $D\in H^{J,a,*}$ satisfies the following conditions:
\begin{align}
&\mbox{$D=\sum_{y\in W^J, y\geqq w}r_yS^{J,a}_y$ with $r_w=1$ and $r_y\in
\bigoplus_{i=0}^{\ell(y)-\ell(w)-1}R_i$}\\
&\mbox{for $y\in W^J$ satisfying $w<y$.}\nonumber\\
&\overline{D}=q^{\ell(w)}D.
\end{align}
Then we have $D=D^{J,a}_w$.
\end{proposition}
\begin{proposition}
[Soergel~\cite{Soergel}]
\label{prop:comparison-Q}
Let $w, y\in W^J$ such that $w\leqq y$.
\begin{itemize}
\item[\rm(i)]
We have $Q^{J,-1}_{w,y}=Q_{w,y}$.
\item[\rm(ii)]
If $W_J$ is a finite group, then we have
\[
Q^{J,q}_{w,y}=\sum_{x\in W_J, ww_J\leqq yx}(-1)^{\ell(x)+\ell(w_J)}Q_{ww_J,yx}.
\]
\end{itemize}
\end{proposition}

\section{Hodge modules}
\label{sec:HM}
\setcounter{equation}{0}
In this section we briefly recall the notation from the theory of Hodge
modules due to M.\ Saito~\cite{S}.

We denote by $\HS$ the category of mixed Hodge structures and by $\HS_k$ the
category of pure Hodge structures with weight $k\in\BZ$.
Let $R$ and $R_k$ be the Grothendieck groups of $\HS$ and $\HS_k$ respectively.
Then we have $R=\bigoplus_{k\in\BZ}R_k$ and $R$ is endowed with a structure
of a commutative ring via the tensor product of mixed Hodge structures.
The identity element of $R$ is given by $[\BQ^H]$, where $\BQ^H$ is the
trivial Hodge structure.
We denote by $R\ni r\mapsto \overline{r}\in R$ the involutive ring
endomorphism induced by the duality functor $\BD:\HS\to\HS^{\op}$.
Here $\HS^{\op}$ denotes the opposite category of $\HS$.
Let $\BQ^H(1)$ and $\BQ^H(-1)$ be the Hodge structure of Tate and its dual
respectively, and set $\BQ^H(\pm n)=\BQ^H(\pm1)^{\otimes n}$ for
$n\in\BZ_{\geq0}$.
We can regard $\BZ[q,q^{-1}]$ as a subring of $R$ by $q^n=[\BQ^H(-n)]$.
Then the condition \eqref{eq:cond:R} is satisfied for this $R$.

Let $Z$ be a finite-dimensional algebraic variety over $\BC$.
There are two conventions for perverse sheaves on $Z$ according to whether
$\BQ_U[\dim U]$ is a perverse sheaf or $\BQ_U[-\codim U]$ is a perverse
sheaf for a closed smooth subvariety $U$ of $Z$.
Correspondingly, we have two conventions for Hodge modules.
When we use the convention so that $\BQ_U[\dim U]$ is a perverse sheaf, we
denote the category of Hodge modules on $Z$ by $\HM(Z)$, and when we use
the other one we denote it by $\tHM(Z)$.
Let $D^b(\HM(Z))$ and $D^b(\tHM(Z))$ denote the bounded derived categories
of $\HM(Z)$ and $\tHM(Z)$ respectively.
Note that $d$ is for dimension and $c\,$ for codimension.
Then the functor $\HM(Z)\to\tHM(Z)$ given by $M\mapsto M[-\dim Z]$ gives the
category equivalences
\[
\HM(Z)\cong\tHM(Z), \qquad D^b(\HM(Z))\cong D^b(\tHM(Z)).
\]
We shall identify $D^b(\HM(Z))$ with $D^b(\tHM(Z))$ via this equivalence,
and then we have
\begin{equation}
\label{eq:HM:equivalence}
\tHM(Z)=\HM(Z))[-\dim Z].
\end{equation}
Although there are no essential differences between $\HM(Z)$ and $\tHM(Z)$,
we have to be careful in extending the theory of Hodge modules to the
infinite-dimensional situation.
In dealing with sheaves supported on finite-dimensional subvarieties
embedded into an infinite-dimensional manifold we have to use $\HM$, while
we need to use $\tHM$
when we treat sheaves supported on finite-codimensional subvariety of an
infinite-dimensional manifold.
In fact what we really need in the sequel is the results for
infinite-dimensional situation; however, we shall only give below a brief
explanation for the finite-dimensional case.
The extension of $\HM$ to the infinite-dimensional situation dealing with
sheaves supported on finite-dimensional subvarieties is easy, and as for the
extension of $\tHM$ to the infinite-dimensional situation dealing with
sheaves supported on finite-codimensional subvarieties we refer the readers
to Kashiwara-Tanisaki~\cite{KT1}.

Let $Z$ be a finite-dimensional algebraic variety over $\BC$.
When $Z$ is smooth, one has a Hodge module $\BQ^H_Z[\dim Z]\in\Ob(\HM(Z))$
corresponding to the perverse sheaf $\BQ_Z[\dim Z]$.
More generally, for a locally closed smooth subvariety $U$ of $Z$  one has a
Hodge module ${}^\pi\BQ^H_U[\dim U]\in\Ob(\HM(Z))$ corresponding to the
perverse sheaf ${}^\pi\BQ_U[\dim U]$.
For $M\in\Ob(D^b(\HM(Z)))$ and $n\in\BZ$ we set $M(n)=M\otimes\BQ^H(n)$.
One has the duality functor
\begin{equation}
\BDD:\HM(Z)\to\HM(Z)^\op, \qquad \BDD:D^b(\HM(Z))\to D^b(\HM(Z))^\op
\end{equation}
satisfying $\BDD\circ\BDD=\Id$,
and we have
\begin{equation}
\BDD({}^\pi\BQ^H_U[\dim U])={}^\pi\BQ^H_U[\dim U](\dim U)
\end{equation}
for a locally closed smooth subvariety $U$ of $Z$.

Let $f:Z\to Z'$ be a morphism of finite-dimensional algebraic varieties.
Then one has the functors:
\begin{align*}
&f^*:D^b(\HM(Z'))\to D^b(\HM(Z)), \qquad f^!:D^b(\HM(Z'))\to D^b(\HM(Z)),\\
&f_*:D^b(\HM(Z))\to D^b(\HM(Z')), \qquad f_!:D^b(\HM(Z))\to D^b(\HM(Z')),
\end{align*}
satisfying
\[
f^*\circ\BDD=\BDD\circ f^!, \qquad f_*\circ\BDD=\BDD\circ f_!.
\]

We define the functors $f^*, f^!, f_*, f_!$ for $D^b(\tHM)$ by identifying
$D^b(\tHM)$ with $D^b(\HM)$.
For $\tHM$ we use the modified duality functor
\begin{equation}
\BDC:\tHM(Z)\to\tHM(Z)^\op, \qquad \BDC:D^b(\HM(Z))\to D^b(\HM(Z))^\op
\end{equation}
given by
\[
\BDC(M)=(\BDD(M))[-2\dim Z](-\dim Z).
\]
It also satisfies $\BDC\circ\BDC=\Id$.
For a locally closed smooth subvariety $U$ of $Z$ we have
${}^\pi\BQ^H_U[-\codim U]\in\Ob(\tHM(Z))$ and
\begin{equation}
\BDC({}^\pi\BQ^H_U[-\codim U])={}^\pi\BQ^H_U[-\codim U](-\codim U).
\end{equation}

When $f:Z\to Z'$ is a proper morphism, we have $f_*=f_!$ and hence
$f_!\circ\BDD=\BDD\circ f_!$.
When $f$ is a smooth morphism, we have $f^!=f^*[2(\dim Z-\dim Z')](\dim
Z-\dim Z')$ and hence $f^*\circ\BDC=\BDC\circ f^*$.

\section{Finite-codimensional Schubert varieties}
\label{sec:FC}
\setcounter{equation}{0}
Let $\Gg$ be a symmetrizable Kac-Moody Lie algebra over $\BC$.
We denote by $W$ its Weyl group and by $S$ the set of simple roots.
Then $(W,S)$ is a Coxeter system.
We shall consider the Hecke algebra $H=H(W)$ over the Grothendieck ring $R$
of the category $\HS$ (see \S~\ref{sec:HM}), and use the notation in
\S~\ref{sec:KL}

Let $X=G/B$ be the flag manifold for $\Gg$ constructed in Kashiwara~\cite{K1}.
Here $B$ is the ``Borel subgroup" corresponding to the standard Borel
subalgebra of $\Gg$.
Then $X$ is a scheme over $\BC$ covered by open subsets isomorphic to
\[
\BA^\infty=\Spec\BC[x_k;k\in\BN]
\]
(unless $\dim\Gg<\infty$).

Let $1_X=eB\in X$ denote the origin of $X$.
For $w\in W$ we have a point $w1_X=wB/B\in X$.
Let $B^-$ be the ``Borel subgroup" opposite to $B$, and set
$X^w=B^-w1_X=B^-wB/B$ for $w\in W$.
Then we have the following result.
\begin{proposition}
[Kashiwara~\cite{K1}]
\label{prop:FC:Schubert}
\begin{itemize}
\item[\rm(i)]
We have $X=\bigsqcup_{w\in W}X^w$.
\item[\rm(ii)]
For $w\in W$, $X^w$ is a locally closed subscheme of $X$ isomorphic to
$\BA^\infty$ $($unless $\dim\Gg<\infty)$ with codimension $\ell(w)$.
\item[\rm(iii)]
For $w\in W$, we have $\overline{X^w}=\bigsqcup\limits_{y\in W, y\geqq w}X^y$.
\end{itemize}
\end{proposition}
We call $X^w$ for $w\in W$ a finite-codimensional Schubert cell, and
$\overline{X^w}$ a finite-codimensional Schubert variety.

Let $J$ be a subset of $S$.
We denote by $Y$ the partial flag manifold corresponding to $J$.
Let $\pi:X\to Y$ be the canonical projection and set $1_Y=\pi(1_X)$.
We have $\pi(w1_X)=1_Y$ for any $w\in W_J$.
For $w\in W^J$ we set $Y^w=B^{-}w1_Y=\pi(X^w)$.
When $W_J$ is a finite group, we have $Y=G/P_J$ and $Y^w=B^-wP_J/P_J$, where
$P_J$ is the ``parabolic subgroup" corresponding to $J$ (we cannot define
$P_J$ as a group scheme unless $W_J$ is a finite group).

Similarly to Proposition~\ref{prop:FC:Schubert} we have the following.
\begin{proposition}
\label{prop:FC:YSchubert}
\begin{itemize}
\item[\rm(i)]
We have $Y=\bigsqcup_{w\in W^J}Y^w$.
\item[\rm(ii)]
For $w\in W^J$, $Y^w$ is a locally closed subscheme of $Y$ isomorphic to
$\BA^\infty$ $($unless $\dim Y<\infty)$ with codimension $\ell(w)$.
\item[\rm(iii)]
For $w\in W^J$, we have $\overline{Y^w}=
\bigsqcup\limits_{y\in W^J, y\geqq w}Y^y$.
\item[\rm(iv)]
For $w\in W^J$, we have $\pi^{-1}(Y^w)=\bigsqcup_{x\in W_J}X^{wx}$.
\end{itemize}
\end{proposition}


We call a subset $\Omega$ of $W^J$ (resp.\ $W$) admissible if it satisfies
\begin{equation}
w, y\in W^J \mbox{(resp.\ $W$)}, w\leqq y, y\in \Omega\Rightarrow w\in \Omega.
\end{equation}
For a finite admissible subset $\Omega$ of $W^J$ we set
$Y^\Omega=\bigcup_{w\in\Omega}Y^w$.
It is a quasi-compact open subset of $Y$.
Let $\tHMB(Y^\Omega)$ be the category of ${B^-}$-equivariant Hodge modules
on $Y^\Omega$ (see Kashiwara-Tanisaki~\cite{KT1} for the equivariant Hodge
modules on infinite-dimensional manifolds), and denote its Grothendieck
group by $K(\tHMB(Y^\Omega))$.
For $w\in W^J$ the Hodge modules $\BQ^H_{Y^w}[-\ell(w)]$ and
${}^\pi\BQ^H_{Y^w}[-\ell(w)]$ are objects of $K(\tHMB(Y^\Omega))$.
Note that $\BQ_{Y^w}[-\ell(w)]$ is a perverse sheaf on $Y$
because $Y^w$ is affine.
Set
\begin{equation}
\tHMB(Y)=\mathop{\varprojlim}\limits_{\Omega}\tHMB(Y^\Omega),\,
K(\tHMB(Y))=\mathop{\varprojlim}\limits_{\Omega}K(\tHMB(Y^\Omega)),
\end{equation}
where $\Omega$ runs through finite admissible subsets of $W^J$.
By the tensor product, $K(\tHMB(Y))$ is endowed with a structure of an
$R$-module.
Then any element of $K(\tHMB(Y))$ is uniquely written as an infinite sum
\[
\sum_{w\in W^J}r_w[\BQ^H_{Y^w}[-\ell(w)]] \mbox{ with } r_w\in R.
\]
Denote by $K(\tHMB(Y))\ni m\mapsto \overline{m}\in K(\tHMB(Y))$ the
involution induced by the duality functor $\BDC$.
Then we have $\overline{rm}=\overline{r}\,\overline{m}$ for any $r\in R$ and
$m\in K(\tHMB(Y))$.

We can similarly define $\tHMB(X)$, $\BQ^H_{X^w}[-\ell(w)]$ and
${}^\pi\BQ^H_{X^w}[-\ell(w)]$  for $w\in W$, $K(\tHMB(X))$, and
$K(\tHMB(X))\ni m\mapsto \overline{m}\in K(\tHMB(X))$ (for $J=\emptyset$).

Let $\pt$ denote the algebraic variety consisting of a single point.
For $w\in W$ (resp.\ $w\in W^J$) we denote by $i_{X,w}:\pt\to X$ (resp.\
$i_{Y,w}:\pt\to Y$) denote the morphism with image $\{w1_X\}$ (resp.\
$\{w1_Y\}$).
We define homomorphisms
\begin{equation}
\Phi:K(\tHMB(X))\to H^*, \qquad \Phi^J:K(\tHMB(Y))\to H^{J,-1,*}
\end{equation}
of $R$-modules by
\begin{align}
&\Phi([M])=\sum_{w\in W}\left(\sum_{k\in\BZ}(-1)^k[H^ki_{X,w}^*(M)]\right)S_w,
\\
&\Phi^J([M])=\sum_{w\in W^J}
\left(\sum_{k\in\BZ}(-1)^k[H^ki_{Y,w}^*(M)]\right)S^{J,-1}_w.
\end{align}

By the definition we have
\begin{align}
\label{eq:Phi:X}
&\Phi([\BQ^H_{X^w}[-\ell(w)]])=(-1)^{\ell(w)}S_w\quad\mbox{ for $w\in W$},\\
\label{eq:Phi:Y}
&\Phi^J([\BQ^H_{Y^w}[-\ell(w)]])=(-1)^{\ell(w)}S^{J,-1}_w\quad\mbox{ for
$w\in W^J$},
\end{align}
and hence $\Phi$ and $\Phi^J$ are isomorphisms of $R$-modules.

The projection $\pi:X\to Y$ induces a homomorphism
\[
\pi^*:K(\tHMB(Y))\to K(\tHMB(X))
\]
of $R$-modules.
\begin{lemma}
\label{lem:FC:Y:dual}
\begin{itemize}
\item[\rm(i)]
The following diagram is commutative.
\[
\begin{CD}
K(\tHMB(Y)) @>{\Phi^J}>>H^{J,-1,*}\\
@V{\pi^*}VV
@VV{{}^t\varphi^{J,-1}}V\\
K(\tHMB(X)) @>>{\Phi}>H^{*}
\end{CD}
\]
\item[\rm(ii)]
$\overline{\pi^*(m)}=\pi^*(\overline{m})$ for any $m\in K(\tHMB(Y))$.
\item[\rm(iii)]
$\overline{\Phi(m)}=\Phi(\overline{m})$ for any $m\in K(\tHMB(X))$.
\item[\rm(iv)]
$\overline{\Phi^J(m)}=\Phi^J(\overline{m})$ for any $m\in K(\tHMB(Y))$.
\end{itemize}
\end{lemma}
\begin{proof}
For $w\in W^J$  we have
$\pi^*(\BQ^H_{Y^{w}})=\BQ^H_{\pi^{-1}Y_w}$, and hence
Proposition \ref{prop:FC:YSchubert} (iv) implies
\[
\pi^*([\BQ^H_{Y^{w}}])=
\sum_{x\in W_J}[\BQ^H_{X^{wx}}].
\]
Thus (i) follows from \eqref{eq:Phi:X}, \eqref{eq:Phi:Y} and \eqref{eq:SJa-sum}

Locally on $X$ the morphism $\pi$ is a projection of the form
$Z\times\BA^\infty \to Z$, and thus $\pi^*\circ\BDC=\BDC\circ\pi^*$.
Hence the statement (ii) holds.

The statement (iii) is already known (see Kashiwara-Tanisaki~\cite{KT1}).

Then the statement (iv) follows from (i), (ii), (iii), \eqref{eq:t-phi-bar}
and the injectivity of ${}^t\varphi^{J,-1}$.
\end{proof}
\begin{theorem}
\label{thm:FC:Y:GG}
Let $w, y\in W^J$ satisfying $w\leqq y$.
Then we have
\[
H^{2k+1}i_{Y,y}^*({}^\pi\BQ^H_{Y^w})=0,\qquad
H^{2k}i_{Y,y}^*({}^\pi\BQ^H_{Y^w})=\BQ^H(-k)^{\oplus Q^{J,-1}_{w,y,k}}
\]
for any $k\in\BZ$.
In particular, we have
\[
\Phi^J([{}^\pi\BQ^H_{Y^w}[-\ell(w)]])=(-1)^{\ell(w)}D^{J,-1}_w.
\]
\end{theorem}
\begin{proof}
Let $w\in W^J$ and set
\[
(-1)^{\ell(w)}\Phi^J([{}^\pi\BQ^H_{Y^w}[-\ell(w)]])
=D
=\sum_{y\in W^J, y\geqq w}r_yS^{J,-1}_y.
\]
By the definition of ${}^\pi\BQ^H_{Y^w}[-\ell(w)]$ we have
\[
\BDC({}^\pi\BQ^H_{Y^w}[-\ell(w)])={}^\pi\BQ^H_{Y^w}[-\ell(w)](-\ell(w)),
\]
and hence we obtain
\begin{equation}
\label{eq:FC:cond1}
\overline{D}=q^{\ell(w)}D
\end{equation}
by Lemma~\ref{lem:FC:Y:dual} (iv).
By the definition of $\Phi^J$ we have
\begin{equation}
\label{eq:FC:cond2}
r_y=\sum_{k\in\BZ}(-1)^k[H^{k}i_{Y,y}^*({}^\pi\BQ^H_{Y^w})],
\end{equation}
and by the definition of ${}^\pi\BQ^H_{Y^w}[-\ell(w)]$ we have
\begin{align}
\label{eq:FC:cond3}
&r_w=1,\\
\label{eq:FC:cond4}
&\mbox{for $y>w$ we have $H^{k}i_{Y,y}^*({}^\pi\BQ^H_{Y^w})=0$ unless}\\
&0\leqq k\leqq (\ell(y)-\ell(w)-1).\nonumber
\end{align}
By the argument similar to Kashiwara-Tanisaki~\cite{KT1} (see also
Kazhdan-Lusztig \cite{KL2}) we have
\begin{equation}
\label{eq:FC:cond5}
[H^{k}i_{Y,y}^*({}^\pi\BQ^H_{Y^w})]\in R_k.
\end{equation}
In particular, we have
\begin{equation}
\label{eq:FC:cond6}
\mbox{for $y>w$ we have $r_y\in\bigoplus_{k=0}^{\ell(y)-\ell(w)-1}R_k$.}
\end{equation}
Thus we obtain $D=D_w^{J,-1}$ by \eqref{eq:FC:cond1}, \eqref{eq:FC:cond3},
\eqref{eq:FC:cond6} and Proposition~\ref{prop:inv-rel-KL}.
Hence $r_y=Q^{J,-1}_{y,w}$.
By \eqref{eq:FC:cond2} and \eqref{eq:FC:cond5} we have
$[H^{2k+1}i_{Y,y}^*({}^\pi\BQ^H_{Y^w})]=0$ and
$[H^{2k}i_{Y,y}^*({}^\pi\BQ^H_{Y^w})]=q^kQ_{w,y,k}$ for any $k\in\BZ$.
The proof is complete.
\end{proof}
By \eqref{eq:SJa-sum-b} and Theorem~\ref{thm:FC:Y:GG} we obtain the following.
\begin{corollary}
\label{cor:FC:Y:mult}
We have
\[
[\BQ^H_{Y^w}[-\ell(w)]]=\sum_{y\geqq
w}P^{J,-1}_{w,y}[{}^\pi\BQ^H_{Y^y}[-\ell(y)]]
\]
in the Grothendieck group $K(\tHMB(Y))$.
In particular,
the coefficient $P^{J,-1}_{w,y,k}$ of
the parabolic Kazhdan-Lusztig polynomial $P^{J,-1}_{w,y}$
is non-negative and equal to
the multiplicity of the irreducible Hodge module
${}^\pi\BQ^H_{Y^y}[-\ell(y)](-k)$ in the
Jordan H\"older series of the Hodge module $\BQ^H_{Y^w}[-\ell(w)]$.
\end{corollary}
\section{Finite-dimensional Schubert varieties}
\label{sec:FD}
\setcounter{equation}{0}
Set
\begin{equation}
X_w=Bw1_X=BwB/B\quad\mbox{for $w\in W$}.
\end{equation}
Then we have the following result.

\begin{proposition}
[Kashiwara-Tanisaki~\cite{KT2}]
\label{prop:FD:Schubert}
Set $\X=\bigcup_{w\in W}X_w$.
Then $\X$ is the flag manifold considered by Kac-Peterson~{\rm\cite{KP}},
Tits~{\rm\cite{T}}, et al.
In particular, we have the following.
\begin{itemize}
\item[\rm(i)]
We have $\X=\bigsqcup_{w\in W}X_w$.
\item[\rm(ii)]
For $w\in W$ $X_w$ is a locally closed subscheme of $X$ isomorphic to
$\BA^{\ell(w)}$.
\item[\rm(iii)]
For $w\in W$ we have $\overline{X}_w=\bigsqcup\limits_{y\in W, y\leqq w}X_y$.
\end{itemize}
\end{proposition}
We call $X_w$ for $w\in W$ a finite-dimensional Schubert cell and
$\overline{X}_w$ a finite-dimensional Schubert variety.
Note that $\X$ is not a scheme but an inductive limit of finite-dimensional
projective schemes (an ind-scheme).

For $w\in W^J$, we set $Y_w=Bw1_Y=\pi(X_w)$.
Similarly to Proposition~\ref{prop:FD:Schubert} we have the following.
\begin{proposition}
\label{prop:FD:YSchubert}
Set $\Y=\bigcup_{w\in W^J}Y_w$.
Then we have the following.
\begin{itemize}
\item[\rm(i)]
We have $\Y=\bigsqcup_{w\in W^J}Y_w$.
\item[\rm(ii)]
For $w\in W^J$, $Y_w$ is a locally closed subscheme of $Y$ isomorphic to
$\BA^{\ell(w)}$.
\item[\rm(iii)]
For $w\in W^J$, we have $\overline{Y}_w
=\bigsqcup\limits_{y\in W^J, y\leqq w}Y_y$.
\item[\rm(iv)]
For $w\in W^J$, we have $\pi^{-1}(Y_w)=\bigsqcup_{x\in W_J}X_{wx}$.
\end{itemize}
\end{proposition}
For a finite admissible subset $\Omega$ of $W^J$ we set
$Y'_\Omega=\bigcup_{w\in\Omega}Y'_w$.
It is a finite dimensional projective scheme.

Let $\HMB(Y'_\Omega)$ be the category of ${B}$-equivariant Hodge modules on
$Y'_\Omega$.
For $w\in W^J$ the Hodge modules $\BQ^H_{Y_w}[\ell(w)]$ and
${}^\pi\BQ^H_{Y_w}[\ell(w)]$ are objects of $\HMB(Y'_\Omega)$.
Note that $\BQ_{Y_w}[\ell(w)]$ is a perverse sheaf because
$Y_w$ is affine.
Set
\begin{equation}
\HMB(Y')=\mathop{\varinjlim}\limits_{\Omega}\HMB(Y'_\Omega),\,
K(\HMB(Y'))=\mathop{\varinjlim}\limits_{\Omega}K(\HMB(Y'_\Omega)),
\end{equation}
where $\Omega$ runs through finite admissible subsets of $W^J$.
By the tensor product $K(\HMB(Y'))$ is endowed with a structure of an
$R$-module.
Then any element of $K(\HMB(Y'))$ is uniquely written as a finite sum
in two ways
\[
\sum_{w\in W^J}r_w[\BQ^H_{Y_w}[\ell(w)]] \ \mbox{ and }
\sum_{w\in W^J}r_w[{}^\pi\BQ^H_{Y_w}[\ell(w)]]
\ \mbox{ with $r_w$, $r'_w\in R$.}
\]
Denote by $K(\HMB(Y'))\ni m\mapsto \overline{m}\in K(\HMB(Y'))$ the
involution of an abelian group induced by the duality functor $\BDD$.
Then we have $\overline{rm}=\overline{r}\,\overline{m}$ for any $r\in R$ and
$m\in K(\HMB(Y'))$.

We can similarly define $\HMB(X')$, $\BQ^H_{X_w}[\ell(w)]$ and
${}^\pi\BQ^H_{X_w}[\ell(w)]$  for $w\in W$, $K(\HMB(X'))$, and
$K(\HMB(X'))\ni m\mapsto \overline{m}\in K(\HMB(X'))$ (for $J=\emptyset$).

For $w\in W$ (resp.\ $w\in W^J$) we denote by $i_{X',w}:\pt\to X'$ (resp.\
$i_{Y',w}:\pt\to Y'$) denote the morphism with image $\{w1_X\}$ (resp.\
$\{w1_Y\}$).
We define homomorphisms
\begin{equation}
\Psi:K(\HMB(X'))\to H, \qquad \Psi^J:K(\HMB(Y'))\to H^{J,q}
\end{equation}
of $R$-modules by
\begin{align}
&\Psi([M])=\sum_{w\in W}\left(\sum_{k\in\BZ}(-1)^k[H^ki_{X',w}^*(M)]\right)T_w,
\\
&\Psi^J([M])=\sum_{w\in W^J}
\left(\sum_{k\in\BZ}(-1)^k[H^ki_{Y',w}^*(M)]\right)T^{J,q}_w.
\end{align}
By the definition we have
\begin{align}
\label{eq:psi:X}
&\Psi([\BQ^H_{X_w}[\ell(w)]])=(-1)^{\ell(w)}T_w\quad\mbox{ for $w\in W$},\\
\label{eq:psi:Y}
&\Psi^J([\BQ^H_{Y_w}[\ell(w)]])=(-1)^{\ell(w)}T^{J,q}_w\quad\mbox{ for $w\in
W^J$},
\end{align}
and hence $\Psi$ and $\Psi^J$ are isomorphisms.

Let $\pi':X'\to Y'$ denote the projection.
Let $\Omega$ be a finite admissible subset of $W$ and set
$\Omega'=\{w\in W^J\,; \,wW_J\cap\Omega\ne\emptyset\}$.
Then $\Omega'$ is a finite admissible subset of $W^J$ and $\pi'$ induces a
surjective projective morphism $X'_\Omega\to Y'_{\Omega'}$.
Hence we can define a homomorphism $\pi'_!:K(HM^B(X'))\to K(HM^B(Y'))$ of
$R$-modules by
\begin{equation}
\pi'_!([M])=\sum_{k\in\BZ}(-1)^k[H^k\pi'_!(M)].
\end{equation}
\begin{lemma}
\label{lem:FD:Y:dual}
\begin{itemize}
\item[\rm(i)]
The following diagram is commutative.
\[
\begin{CD}
K(\HMB(X')) @>{\Psi}>>H\\
@V{\pi'_!}VV
@VV{\varphi^{J,q}}V\\
K(\HMB(Y')) @>>{\Psi^J}>H^{J,q}
\end{CD}
\]
\item[\rm(ii)]
$\overline{\pi'_!(m)}=\pi'_!(\overline{m})$ for any $m\in K(\HMB(X'))$.
\item[\rm(iii)]
$\overline{\Psi(m)}=\Psi(\overline{m})$ for any $m\in K(\HMB(X'))$.
\item[\rm(iv)]
$\overline{\Psi^J(m)}=\Psi^J(\overline{m})$ for any $m\in K(\HMB(Y'))$.
\end{itemize}
\end{lemma}
\begin{proof}
Let $w\in W^J$ and $x\in W_J$.
Since $X_{wx}\to Y_w$ is an $\BA^{\ell(x)}$-bundle, we have
$\pi'_!(\BQ^H_{X_{wx}})=\BQ^H_{Y_w}[-2\ell(x)](-\ell(x))$, and hence
\[
\pi'_!([\BQ^H_{X_{wx}}[\ell(wx)]])=(-q)^{\ell(x)}[\BQ^H_{Y_w}[\ell(w)]].
\]
Thus (i) follows from \eqref{eq:psi:X}, \eqref{eq:psi:Y} and \eqref{eq:phi-T}.

The statement (ii) follows from the fact that $\pi'$ is an inductive limit
of projective morphisms and hence $\pi'_!$ commutes with the duality functor
$\BDD$.

The statement (iii) is proved similarly to Kashiwara-Tanisaki~\cite{KT1}, and
we omit the details (see also Kazhdan-Lusztig~\cite{KL2}).
Then the statement (iv) follows from (i), (ii), (iii), \eqref{eq:phi,j} and
surjectivity of $\varphi^{J,q}$.
\end{proof}
\begin{theorem}
\label{thm:FD:Y:GG}
Let  $w, y\in W^J$ such that $w\geqq y$.
Then we have
\[
H^{2k+1}i_{Y',y}^*({}^\pi\BQ^H_{Y_w})=0,\qquad
H^{2k}i_{Y',y}^*({}^\pi\BQ^H_{Y_w})=\BQ^H(-k)^{\oplus P^{J,q}_{y,w,k}}
\]
for any $k\in\BZ$.
In particular, we have
\[
\Psi^J([{}^\pi\BQ^H_{Y_w}[\ell(w)]])=(-1)^{\ell(w)}C^{J,q}_w.
\]
\end{theorem}
\begin{proof}
Let $w\in W^J$ and set
\[
(-1)^{\ell(w)}\Psi^J([{}^\pi\BQ^H_{Y_w}[\ell(w)]])
=C
=\sum_{y\in W^J, y\leqq w}r_yT^{J,q}.
\]
By the definition of ${}^\pi\BQ^H_{Y_w}[\ell(w)]$ we have
$\BDD({}^\pi\BQ^H_{Y_w}[\ell(w)])={}^\pi\BQ^H_{Y_w}[\ell(w)](\ell(w))$.
Hence we obtain
\begin{equation}
\label{eq:FD:cond1}
\overline{C}=q^{-\ell(w)}C
\end{equation}
by Lemma~\ref{lem:FD:Y:dual} (iv).
By the definition of $\Psi^J$ we have
\begin{equation}
\label{eq:FD:cond2}
r_y=\sum_{k\in\BZ}(-1)^k[H^{k}i_{Y',y}^*({}^\pi\BQ^H_{Y_w})],
\end{equation}
and by the definition of ${}^\pi\BQ^H_{Y_w}[\ell(w)]$ we have
\begin{align}
\label{eq:FD:cond3}
&r_w=1,\\
\label{eq:FD:cond4}
&\mbox{for $y<w$ we have $H^{k}i_{Y',y}^*({}^\pi\BQ^H_{Y_w})=0$ unless}\\
&0\leqq k\leqq (\ell(w)-\ell(y)-1).\nonumber
\end{align}
Moreover, by the argument similar to Kazhdan-Lusztig~\cite{KL2} and
Kashiwara-Tanisaki~\cite{KT1} we have
\begin{equation}
\label{eq:FD:cond5}
[H^{k}i_{Y',y}^*({}^\pi\BQ^H_{Y_w})]\in R_k.
\end{equation}
In particular, we have
\begin{equation}
\label{eq:FD:cond6}
\mbox{for $y<w$ we have $r_y\in\bigoplus_{k=0}^{\ell(w)-\ell(y)-1}R_k$.}
\end{equation}
Thus we obtain $C=C_w^{J,q}$ by \eqref{eq:FD:cond1}, \eqref{eq:FD:cond3},
\eqref{eq:FD:cond6} and Proposition~\ref{prop:rel-KL}.
Hence $r_y=P^{J,q}_{y,w}$.
By \eqref{eq:FD:cond2} and \eqref{eq:FD:cond5} we have
$[H^{2k+1}i_{Y',y}^*({}^\pi\BQ^H_{Y_w})]=0$ and
$[H^{2k}i_{Y',y}^*({}^\pi\BQ^H_{Y_w})]=q^kP_{y,w,k}$ for any $k\in\BZ$.
The proof is complete.
\end{proof}
We note that a result closely related to Theorem~\ref{thm:FD:Y:GG} above is
already given in Deodhar~\cite{D}.

By \eqref{eq:TJa-sum} and Theorem~\ref{thm:FD:Y:GG} we obtain the following.
\begin{corollary}
\label{cor:FD:Y:mult}
We have
\[
[\BQ^H_{Y'_w}[\ell(w)]]=\sum_{y\leqq
w}Q^{J,q}_{y,w}[{}^\pi\BQ^H_{Y'_y}[\ell(y)]]
\]
in $K(\HMB(Y'))$.
In particular,
the coefficient $Q^{J,q}_{y,w,k}$
of the inverse parabolic Kazhdan-Lusztig polynomial
$Q^{J,q}_{y,w}$ is non-negative
and equal to the multiplicity of the irreducible Hodge module
${}^\pi\BQ^H_{Y'_y}[\ell(y)](-k)$ in the
Jordan H\"older series of the Hodge module $\BQ^H_{Y'_w}[\ell(w)]$.
\end{corollary}

\bibliographystyle{unsrt}
\def\same{\,$\raise3pt\hbox to 25pt{\hrulefill}\,$}

\end{document}